\documentclass{amsart}
\usepackage{amssymb}
\usepackage{amsthm}
\usepackage{hyperref}
\usepackage{amsmath}
\usepackage{amssymb}
\usepackage{enumerate}
\usepackage{graphicx}
\usepackage{amsthm}
\usepackage{amsopn}
\usepackage{amsfonts}
\usepackage{upgreek}
\usepackage{amsfonts}
\usepackage{amssymb}
\usepackage{empheq}
\usepackage{caption}
\numberwithin{equation}{section}

\newtheorem{theorem}{Theorem}
\newtheorem{corollary}{Corollary}
\newtheorem{example}{Example}%
\newtheorem{remark}{Remark}%
\newtheorem{lemma}{Lemma}%
\newtheorem{definition}{Definition}%

\begin{document}

\title[On oscillations of nonautonomous linear IDEPCAG]{On oscillations of nonautonomous linear impulsive differential equations with general piecewise constant deviating arguments}
\author{Ricardo Torres}
\address{Instituto de Ciencias Físicas y Matemáticas, Facultad de Ciencias, Universidad Austral de Chile\\
Campus Isla Teja s/n, Valdivia, Chile.}
\address{\noindent Instituto de Ciencias, Universidad Nacional de General Sarmiento\\
Los Polvorines, Buenos Aires, Argentina}
\curraddr{}
\email{ricardo.torres@uach.cl}
\thanks{}



\date{}

\dedicatory{Dedicated to the memory of Prof. Istv\'an Gy\H ori and Prof. Nicol\'as Yus Su\'arez.}

\begin{abstract}
We provide sufficient conditions that ensure oscillations and nonoscillations for nonautonomous impulsive differential equations with piecewise constant arguments of a generalized type. We cover several cases of differential equations with deviated arguments investigated before as particular cases.
\end{abstract}
\maketitle
\section{Introduction}
In many cases, the functions used for modeling the natural phenomena with differential equations may have to be discontinuous, such as a piecewise constant, or have the impulsive effect present to reflect the phenomena's characteristics to be modeled. The reader can find examples of this kind of phenomenon in the works of S. Busenberg and K. Cooke (\cite{COOKE1}) and L. Dai and M.C. Singh (\cite{DAI_SINGH}). In the first one, the authors modeled diseases of vertical transmission. In the second one, the authors studied the oscillatory motion of spring-mass systems subject to piecewise constant forces of the form $Ax([t])$ or $A\cos([t]).$ Indeed, one such system is the famous Geneva wheel, widely used in watches and instruments. The motion of this mechanism is piecewise continuous, and it can be modeled by $$mx''(t)+kx_1=A sin\left(\omega \left[\dfrac{t}{T}\right]\right),$$
where $[\cdot]$ is the greatest integer function. 
 (See \cite{DAI}).\\

As a generalization of the work done by A. Myshkis in \cite{203} (1977) concerning differential equations with deviating arguments, M. Akhmet introduced a very interesting class of differential equations of the form
\begin{equation}
z^{\prime}(t)=f(t,z(t),z(\gamma(t))),\label{depcag_eq}
\end{equation}
where $\gamma(t)$ is a \emph{piecewise constant argument of generalized type}.
Such $\gamma$ is defined as follows: let $\left(t_{n}\right)_{n\in\mathbb{Z}}$ and $\left(\zeta_{n}\right)_{n\in\mathbb{Z}}$ such that $t_{n}<t_{n+1}\, , \forall n\in\mathbb{Z}$ with $\displaystyle{\lim_{n\rightarrow\infty}t_{n}=\infty}$, $\displaystyle{\lim_{n\rightarrow -\infty}t_{n}=-\infty}$ and $\zeta_{n}\in[t_{n},t_{n+1}],$
then  $\gamma(t)=\zeta_{n},$ if $t\in I_{n}=\left[t_{n},t_{n+1}\right).$
I.e., $\gamma(t)$ is a locally constant function. An elementary example of such functions is 
$\gamma(t) =[t],$ where $[\cdot]$ is the greatest integer function, which is constant in every interval $[n,n+1[$ with  $n\in \mathbb{Z}$.\\
If a piecewise constant argument is used, the interval $I_n$ is decomposed into an advanced and delayed subintervals $I_{n}=I_{n}^{+}\bigcup I_{n}^{-}$, where
$I_{n}^{+}=[t_{n},\zeta_{n}]$ and $I_{n}^{-}=[\zeta_{n},t_{n+1}].$
This class of differential equations is known as \emph{Differential Equations with Piecewise Constant Argument of Generalized Type} (\emph{DEPCAG}). They have remarkable properties, as the solutions are continuous functions, even when $\gamma$ is discontinuous. If we assume continuity of the solutions of \eqref{depcag_eq}, integrating from $t_n$ to $t_{n+1}$, we define a difference equation, so this type of differential equation has hybrid dynamics (see \cite{AK2, P2011, Wi93}).\\
For example, in \cite{Torres_proyecciones}, the author introduced the piecewise constant argument $\gamma(t) =\left[ \frac{t}{h}\right]h+\alpha h $ with $h>0$ and $0\leq \alpha\leq 1.$ Then, we can see that
\begin{eqnarray*}
\left[ \frac{t}{h}\right]h+\alpha h=(n+\alpha)h, \text{ when } t\in I_n=[nh,\left(n+1\right)h).
\end{eqnarray*}
Its delayed and advanced intervals can be concluded easily:  $\gamma (t) -t \geq 0\text{ }\Leftrightarrow \text{ }t\leq (n+\alpha)h$ and
$\gamma (t) -t \leq 0\text{ }\Leftrightarrow \text{ }t\geq (n+\alpha)h$. Hence, we have
\begin{equation*}
I_{n}^{+}=[nh,(n+\alpha)h),\quad I_{n}^{-}=[(n+\alpha)h,\left(n+1\right)h).
\end{equation*}
This piecewise constant argument has been used (case $\alpha=0$) in the approximation of solutions of ordinary and delayed differential equations (see also \cite{TORRES3, TORRES5} ).\\

 If an impulsive condition is considered at instants $\{t_n\}_{n\in\mathbb{Z}}$, we define the 
\emph{Impulsive differential equations with piecewise constant argument of generalized type} (\emph{IDEPCAG}),
\begin{eqnarray}
&&z^{\prime}(t)=f(t,z(t),z(\gamma(t))),\qquad \qquad \qquad \,\quad t\neq t_{n}, \nonumber \\
&&\Delta z(t_{n}):=z(t_{n})-z(t_{n}^{-})=P_{n}(z(t_{n}^{-})),\qquad t=t_{n},\quad n\in\mathbb{N} \label{idepcag_gral},
\end{eqnarray}
where $z(t_n^-)=\displaystyle{\lim_{t\to t_n^-}z(t),}$ and $P_n$ is the impulsive operator (see \cite{AK3,Samoilenko}).\\
If the differential equation explicitly shows the piecewise constant argument used, we will call it DEPCA (IDEPCA if it has impulses).\\

\begin{definition}
A function $x(t)$ defined on $[\tau, \infty)$ is said to be oscillatory if there exist two real-valued sequences $(a_n),\,(b_n)\subset [\tau,\infty)$ such that $a_n\to\infty,\,b_n\to\infty$ as $n\to\infty$ and 
$x(a_n)\leq 0\leq x(b_n),\,\, \forall n\geq M,$ where $M$ is sufficiently large. I.e., 
if it is neither eventually positive nor eventually negative. Otherwise, it is called nonoscillatory. 
\end{definition}
\begin{figure}[h!]
\centering
\includegraphics[scale=0.15]{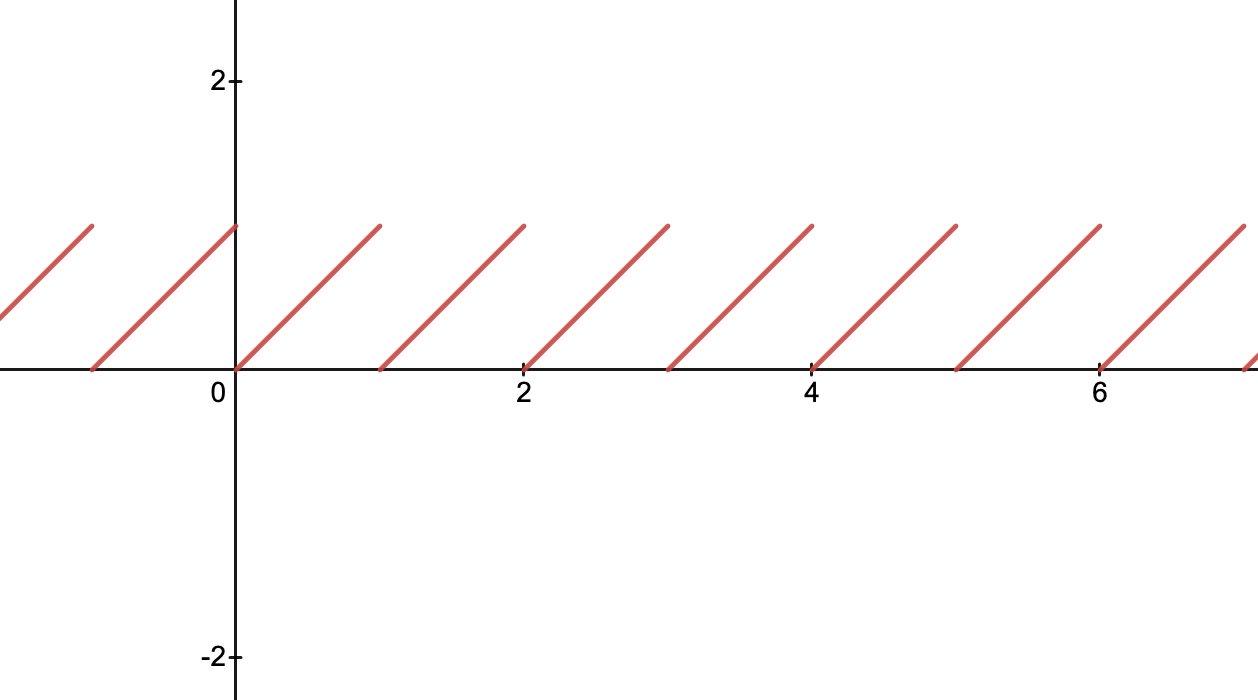}
\caption{$f(t)=t-[t]$: an example of an oscillatory piecewise continuous function with infinite zeros ($f(n)= 0, \,\forall n\in\mathbb{Z}$).}
\end{figure}

In the piecewise continuous frame, a function $x(t)$ can be oscillatory even if $x(t)\neq 0,\, \forall t\in [\tau,\infty)$. \textit{This is the context of this work, so all the definition of oscillation for the continuous frame are no longer valid (see \cite{karakoc_berek})}.
\begin{figure}[h!]
\centering
\includegraphics[scale=0.15]{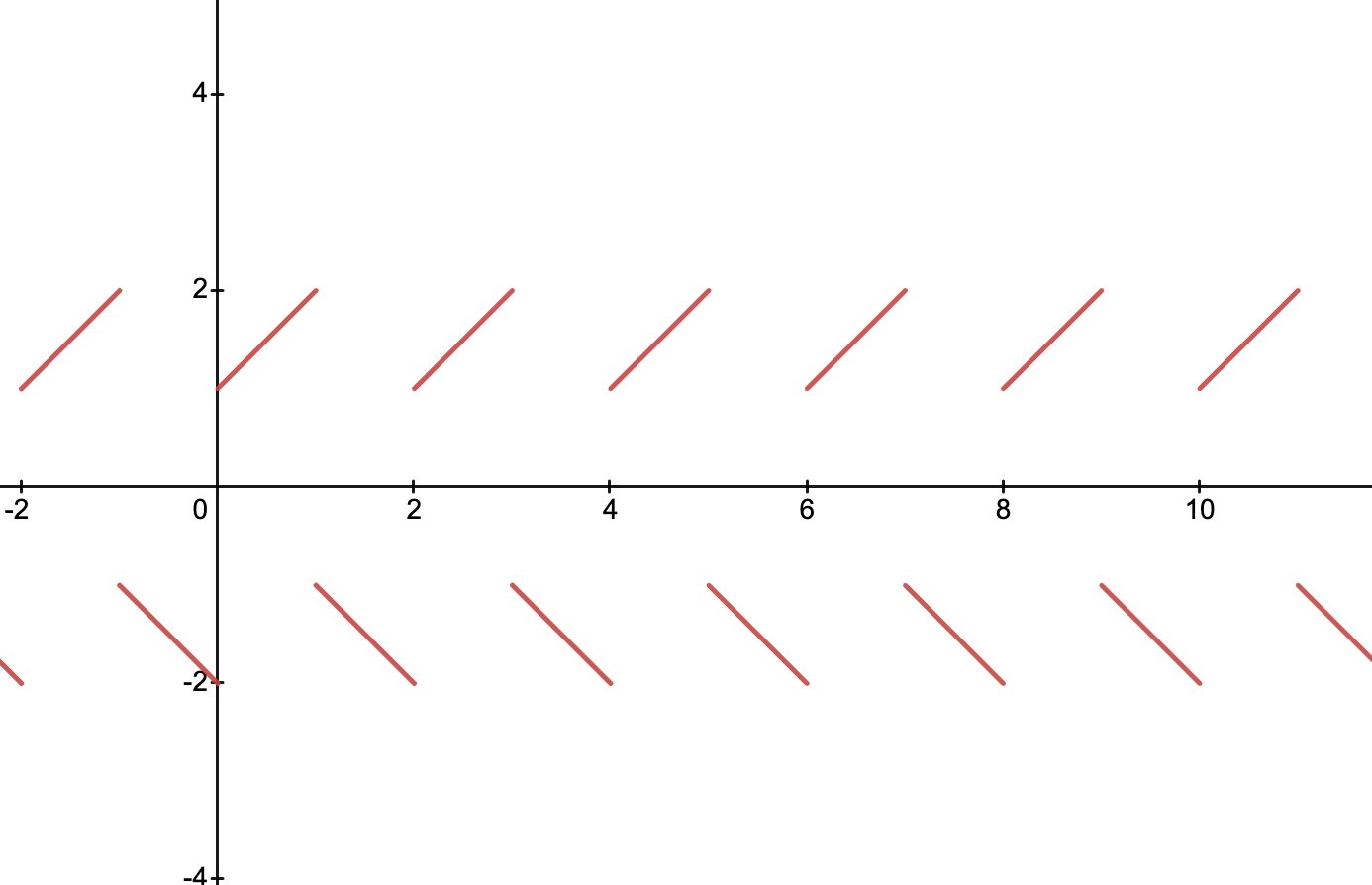}
\caption{$f(t)=(-1)^{[t]}(t+1-[t])$. An oscillatory piecewise continuous function with $f(t)\neq 0, \,\forall t\in\mathbb{R}$.}
\end{figure}

\begin{definition}
A nonzero solution $y(n)$ of a difference equation will be \it{oscillatory} if for every positive integer $M>0$, there exists $n\geq M$ such that 
$$y(n)y(n+1)\leq 0.$$ 
Otherwise, $y(n)$ will be called \it{nonoscillatory} (see \cite{ELAIDY}).
\end{definition}
\section{Aim of the work}
In \cite{WiAf88}, \textit{J.Wiener} and \textit{A.R. Aftabizadeh}
studied the following DEPCA 
\begin{equation}
\begin{tabular}{ll}
$x^{\prime }(t)=a(t)x(t)+a_0(t)x\left(m\left[\frac{t+k}{m}\right]\right),$ & $x(0)=C_0,$ 
\end{tabular}
\label{idepca_afta_wiener}
\end{equation}
where $a(t),a_0(t)$ are continuous functions on $[0,\infty).$ They conclude the following oscillatory criterion:
\begin{theorem}\label{afta_wiener}
if either of the conditions
\begin{eqnarray*}
&&\limsup_{n\to\infty} \displaystyle{\int_{mn-k}^{mn}a_0(t)\exp{\left(\int_{t}^{mn} a(s)ds\right)}dt>1},\\
&&\liminf_{n\to\infty} \displaystyle{\int_{mn}^{m(n+1)-k}a_0(t)\exp{\left(\int_{t}^{mn} a(s)ds\right)}dt<-1},
\end{eqnarray*}    
  holds true, then every solution of \eqref{afta_wiener} is oscillatory.
\end{theorem}
Inspired by the last work, we will get sufficient conditions for the oscillation of the solutions of the scalar linear nonautonomous \emph{IDEPCAG} 
\begin{equation}
\begin{tabular}{ll}
$x^{\prime }(t)=a(t)x(t)+b(t)x(\gamma (t)),$ & $t\neq t_{k}$ \\ 
$\Delta x|_{t=t_{k}}=c_{k}x(t_{k}^{-}),$ & $t=t_{k}$%
\end{tabular}
\label{sistema_idepcag_general_abstract}
\end{equation}
where $a(t),b(t)$ are continuous real-valued locally integrable functions, $c_k\neq -1,\,\forall k\in\mathbb{Z}$ and $\gamma$ is a piecewise constant argument of generalized type.\\
We will improve the results given in \cite{WiAf88} considering any piecewise constant argument $\gamma(t)$ such that $\gamma(t)=\zeta_{k}$ if $t\in[t_k,t_{k+1})$ with $t_k\leq \zeta_k\leq t_{k+1}.$ and the impulsive effect.\\

Our work is organized as follows: First, we provide some preliminary information necessary to understand the IDEPCAG frame better. Next, we study oscillations for the scalar linear IDEPCAG and provide some criteria for oscillatory solutions (extending the classical ones) for the autonomous and nonautonomous cases. Finally, we give some examples showing the effectiveness of our results. 

\section{Some recent results of oscillations on IDEPCAG}
In  \cite{Kuo_pinto_2013} (2013), K-S. Chiu and M. Pinto studied the scalar DEPCAG
$$y'(t)=a(t)y(t)+b(t)y(\gamma(t)),\quad y(\tau)=y_0,$$
where $a(t),b(t)$ are real continuous functions defined on $\mathbb{R}.$ The authors investigated oscillation and periodic solutions. This excellent work seems to be the first to consider oscillations and periodic solutions with a generalized piecewise constant argument. One of the main hypotheses used by the authors is the continuity of the solutions. Still, it has no results for the IDEPCA or IDEPCAG cases.\\

In \cite{karakoc_berek} (2018), F. Karakoc,  A. Unal, and H.  Bereketoglu studied the following nonlinear IDEPCA:
$$
\begin{cases}\label{berek_1}
x'(t)=-a(t)x(t)-x([t-1])f(y([t])+h_1(x[t]), & t\neq n\in\mathbb{Z}^+   \\
y'(t)=-b(t)x(t)-y([t-1])g(y([t])+h_2(y[t]), & 
\end{cases}
$$
$$
\begin{cases}\label{berek_2}
x(n)=(1+c_n)x(n^-), & t= n\in\mathbb{Z}^+   \\
y(n)=(1+d_n)y(n^-), & 
\end{cases}
$$
with initial conditions $x(-1)=x_{-1},\,x(0)=x_0, y(-1)=y_{-1},\,y(0)=y_0,$ where $a,b:[0,\infty)\to\mathbb{R}$ are continuous functions, $f,g,h_1,h_2\in C(\mathbb{R},\mathbb{R}),$ $c_n$ and $d_n$ are sequences of real numbers such that $c_n\neq 1$, $d_n\neq 1,\,\forall n\geq 1.$ The authors studied the existence and uniqueness of solutions, and they obtained some sufficient conditions for the oscillation of the solution. They stated the following sufficient conditions of oscillation of the solutions:
\begin{theorem}\label{Theorem_berek_karakoc}
    Assume that there exist $M_1>0$ and $M_2>0$ such that $f(t)\geq M_1$ and $g(t)\geq M_2\,\, ,\forall t\in\mathbb{R},$ $t\,h_1(t),t\,h_2(t)<0$ for $t\neq 0$ and $c_n,\,d_n<1$ for all $n\in\mathbb{Z}^+.$ If the following conditions hold, 
\begin{eqnarray*}
    \limsup_{n\to\infty} (1-c_n)\displaystyle{\int_n^{n+1}\exp{\left(\int_{n-1}^s a(u)du\right)}ds>\frac{1}{M_1},}\\
    \limsup_{n\to\infty} (1-d_n)\displaystyle{\int_n^{n+1}\exp{\left(\int_{n-1}^s b(u)du\right)}ds>\frac{1}{M_2}.}
\end{eqnarray*}    
    then all solutions of \eqref{berek_1}-\eqref{berek_2} are oscillatory.
\end{theorem}

In  \cite{Kuo_pinto_2013} (2013), K-S. Chiu and M. Pinto studied the scalar DEPCAG
$$y'(t)=a(t)y(t)+b(t)y(\gamma(t)),\quad y(\tau)=y_0,$$
where $a(t),b(t)$ are real continuous functions defined on $\mathbb{R}.$ The authors investigated oscillation and periodic solutions. This excellent work seems to be the first to consider oscillations and periodic solutions with a generalized piecewise constant argument. One of the main hypotheses used by the authors is the continuity of the solutions. Still, it has no results for the IDEPCA or IDEPCAG cases.\\

In \cite{Kuo2023} (2023), K-S. Chiu and I. Berna considered the following impulsive differential equation with a piecewise constant argument
\begin{align}
y^{\prime}(t)&=a(t)y(t)+b(t)y\left(p\left[\frac{t+l}{p}\right]\right),\quad y(\tau)=c_0,\qquad t\neq kp-l \nonumber \\
\Delta y(kp-l)&=d_k y({kp-l}^{-}),\qquad t=kp-l,\quad k\in\mathbb{Z} \label{kuo_homogenea},
\end{align}

\noindent where $a(t)\neq 0$ and $b(t)$ are real-valued continuous functions, $p<l$ and $d_k\neq -1$, $\forall k\in\mathbb{Z}.$ 
The authors obtained an oscillation criterion for solutions of \eqref{kuo_homogenea}.\\

Finally, in \cite{Zhaolin} (2023), Z. Yan and J. Gao studied oscillation and non-oscillation of Runge–Kutta methods for the following linear mixed-type impulsive differential equations with piecewise constant arguments
\begin{align}
&y'(t)=ay(t)+by([t-1])+cy([t])+dy([t+1]), \quad t\neq k,\nonumber \\
&\Delta y(k)=q y(k^{-}),\qquad t=k,\quad k\in\mathbb{N} \label{zhaolin},\\
&y(0)=y_0, y(1)=y_1,\nonumber
\end{align}
where $a,b,c,d,q, y_0,y_1\in\mathbb{R}$
and $q\neq 1$.
They obtained conditions for oscillation and non-oscillation of the numerical solutions. Also, the authors obtained conditions under which Runge–Kutta methods can preserve the oscillation and non-oscillation.

\subsection{Why study oscillations on IDEPCAG?}
\begin{example}\label{ejemplo_berek_avanzado}
Let the following scalar linear DEPCA
\begin{equation}
    x'(t)=a(t)(x(t)-x([t+1]),\quad x(\tau)=x_0,\label{berek_depca}
\end{equation}
and the scalar linear IDEPCA
\begin{equation}
    \begin{tabular}{ll}
$z'(t)=a(t)\left(z(t)-z([t+1])\right),$ & $t\neq k$ \\ 
$z(k)=c_k z(k^{-}),$ & $t=k,\quad k\in\mathbb{Z},$ 
\end{tabular}
\label{berek_sistema_avanzado}
\end{equation}
where $a(t)$ is a continuous locally integrable function. As $\gamma(t)=[t+1],$ we have $t_{k}=k,$ 
$\zeta_k=k+1$ if $t\in[k,k+1],\,k\in\mathbb{Z}$.\\
It is very easy to see, by Theorem \ref{TEO_FORMULA_var_PAram}  (considering $c_k=0,\,\forall k\in\mathbb{Z}$), that the solution of \eqref{berek_depca} is $x(t)=x_0,\,\, \forall t\geq \tau.$
I.e., all the solutions are constant. Hence, \eqref{berek_depca} has no nontrivial oscillatory solutions. \\
On the other hand, the solution of \eqref{berek_sistema_avanzado}  is
$$z(t)=\left(\prod_{j=k(\tau)+1}^{k(t)}c_j\right)z(\tau),\quad t\geq \tau,$$
where $k(t)=k$ is the only integer such that $t\in[k,k+1]$ (see Example 4). Hence, if $c_k<0,\,\forall k\in\mathbb{Z}$
all the solutions are oscillatory. This example shows the role of the presence of the impulsive effect. In fact, our main result applies to DEPCAG and IDEPCAG cases.
\begin{figure}[h!]
\centering
\includegraphics[scale=0.2]{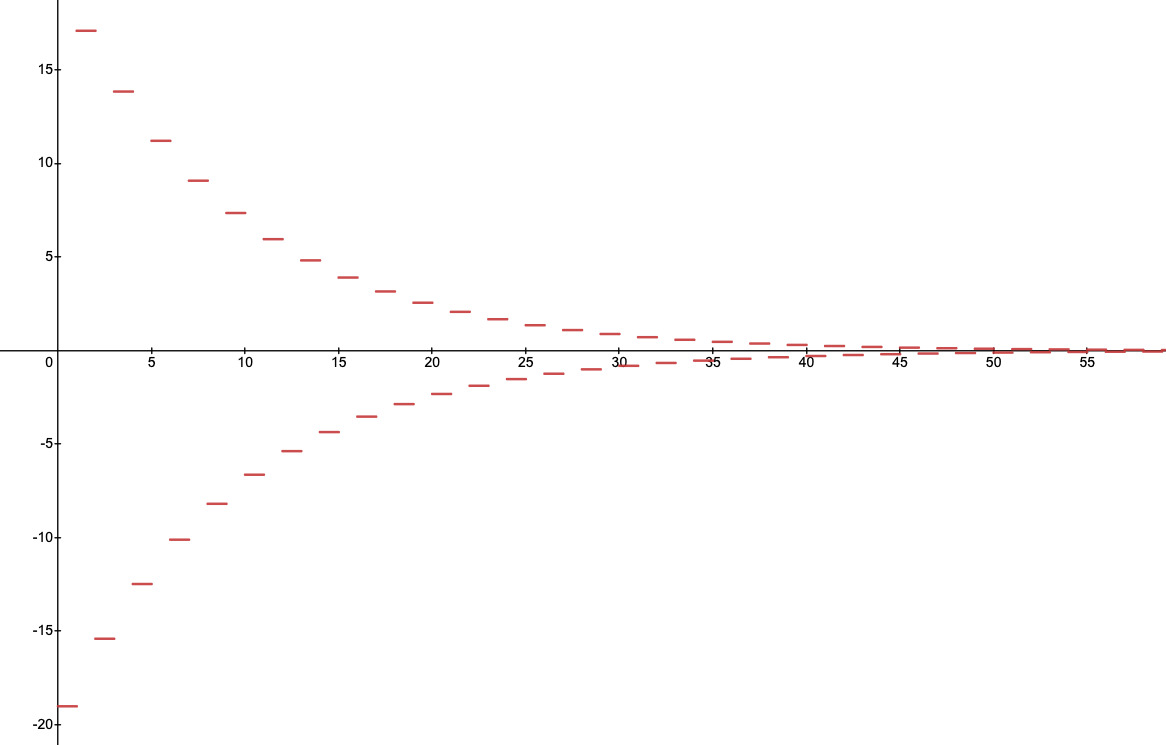}
\caption{Solution of \eqref{berek_sistema_avanzado} with $c_k=-\frac{9}{10}$ and $z(0)=-19.$}
\end{figure}
\end{example}
\section{Some Preliminary on IDEPCAG theory}
Let $\mathcal{PC}(X, Y)$ be the set of all functions $r: X\to Y$ which are continuous for $t\neq t_k$ and continuous from the left with discontinuities of the first kind at $t=t_k$. Similarly, let $\mathcal{PC}^{1}(X,Y)$ the set of functions $s:X\to Y$ such that $s'\in \mathcal{PC}(X,Y).$
\begin{definition}[DEPCAG solution]
A continuous function $x(t)$ is a solution of \eqref{depcag_eq} if:
\begin{itemize}
\item[(i)] $x'(t)$ exists at each point $t\in \mathbb{R}$ with the possible exception at the times $t_{k}$, $k\in \mathbb{Z}$, where the one side derivative exists.
\item[(ii)] $x(t)$ satisfies \eqref{depcag_eq} on the intervals of the form $(t_{k},t_{k+1})$, and it holds for the right
derivative of $x(t)$ at $t_{k}$.
\end{itemize}
\end{definition}
\begin{definition}[IDEPCAG solution]
A piecewise continuous function $y(t)$ is a solution of \eqref{idepcag_gral} if:
\begin{itemize}
\item[(i)] $y(t)$ is continuous on $I_k=[t_k,t_{k+1})$ with first kind discontinuities at $t_k, \,k\in \mathbb{Z}$, 
where $y'(t)$ exists at each $t\in \mathbb{R}$ with the possible exception at the times $t_{k}$, where lateral derivatives exist (i.e. $y(t)\in \mathcal{PC}^{1}([t_k,t_{k+1}),\mathbb{R})$).
\item[(ii)] The ordinary differential equation
$$y^{\prime}(t)=f(t,y(t),y(\zeta_{k}))$$ holds on every interval $I_{k}$, where $\gamma(t)=\zeta_k$.
\item[(iii)]
For $t=t_{k}$, the impulsive condition 
$$\Delta y(t_{k})=y(t_{k})-y(t_{k}^{-})=P_{k}(y(t_{k}^{-}))$$
holds. I.e., $y(t_{k})=y(t_{k}^{-})+P_{k}(y(t_{k}^{-}))$, where $y(t_{k}^{-})$ denotes the left-hand limit
of the function $y$ at $t_k$.
\end{itemize}
\end{definition}
Let the following IDEPCAG:
\begin{equation}
\begin{tabular}{ll}
$x^{\prime }(t)=f(t,x(t),x(\gamma (t))),$ & $t\neq t_{k}$ \\ 
$x(t_{k})-x\left(t_{k}^{-}\right) =P_{k}(x(t_{k}^{-})),$ & $t=t_{k},$ \\ 
$x(\tau )=x_{0},$ & 
\end{tabular}
\label{SISTEMA_GENERICO_IDEPCAG}
\end{equation}
where  $f\in C([\tau ,\infty )\times \mathbb{R}\times \mathbb{R},\mathbb{R}),$  $P_{k}\in C(\left\{t_{k}\right\} ,\mathbb{R})$ and $\left(\tau ,x_{0}\right) \in \mathbb{R}^2$.\\
Also, let the following hypothesis hold:
\begin{itemize}
\item[\textbf{(H1)}] Let $\eta_{1},\eta_{2}:\mathbb{R}\rightarrow [0,\infty )$ locally integrable functions 
and $\omega_{k}\in \mathbb{R}^+$, $\forall k\in \mathbb{Z}$; such that
\begin{eqnarray*}
\left\vert f(t,x_{1},y_{1})-f(t,x_{2},y_{2})\right\vert &\leq &\eta_{1}(t) \left\vert x_{1}-x_{2}\right\vert 
+\eta_{2}(t) \left\vert y_{1}-y_{2}\right\vert ,\\
\left\vert P_{k}(x_{1}(t_{k}^{-}))-P_{k}(x_{2}(t_{k}^{-}))\right\vert 
&\leq &\omega_{k}\left\vert x_{1}\left( t_{k}^{-}\right) -x_{2}\left(t_{k}^{-}\right) \right\vert .
\end{eqnarray*}
\item[\textbf{(H2)}] Assume that
\begin{equation*}
\overline{\nu }=\sup_{k\in \mathbb{Z}}\left(\int_{t_{k}}^{t_{k+1}}\left(\eta_{1}(s)+\eta_{2}(s)\right)ds\right) <1.
\end{equation*}
\end{itemize}

In the following, we mention some valuable results for the rest of the work: 
\subsection{\bf{IDEPCAG Gronwall-Bellman type inequality}}
\begin{lemma}{(\cite{Torres1},\cite{CTP2019} Lemma 4.3)}
\label{LEMA_GRONWALL_PINTO_2}Let $I$ an interval and $u,\eta_{1},\eta_{2}:I\rightarrow [0,\infty )$ such that $u$ is continuous (with possible exception at $\{t_k\}_{k\in\mathbb{N}}$), $\eta_{1},\eta_{2}$ are continuous and locally integrable functions, $\eta:\left\{ t_{k}\right\} \rightarrow [0,\infty )$ and $\gamma (t) $ a piecewise constant argument of generalized type such that $\gamma (t)=\zeta_{k}$,  
$\forall t\in I_{k}=[t_{k},t_{k+1})$ with $t_{k}\leq \zeta_{k}\leq t_{k+1}$ $\forall k\in \mathbb{N}.$ 
Assume that $\forall t\geq \tau $  
\begin{eqnarray*}
&u(t)\leq& u(\tau) +\int_{\tau}^{t}\left(\eta_{1}(s)u(s)+\eta_{2}(s)u(\gamma(s))\right) ds+\sum_{\tau \leq t_{k}<t}\eta(t_{k})u(t_{k}^{-}),\\
&\widehat{\vartheta }_{k}&=\int_{t_{k}}^{\zeta_{k}}\left(\eta_{1}(s)+\eta_{2}(s)\right) ds\leq \widehat{\vartheta }
:=\sup_{k\in \mathbb{N}}\widehat{\vartheta }_{k}<1  \label{CONDICION_INVERTIBILIDAD_GRONWALL_2}
\end{eqnarray*}
hold. Then, for $t\geq \tau $, we have
\begin{align*}
u(t)\leq &\left(\prod_{\tau \leq t_{k}<t}\left(1+\eta(t_{k})\right) \right) 
\exp \left(\int_{\tau}^{t}\left(\eta_{1}(s)+\frac{\eta_{2}(s)}{1-\widehat{\vartheta}}\right) ds\right) u(\tau),\\
u(\zeta_k)\leq &(1-\vartheta)u(t_k)\\
u(\gamma(t))\leq &(1-\vartheta)^{-1}\left(\prod_{\tau \leq t_{k}<t}\left(1+\eta_3 (t_{j})\right) \right) 
\exp \left(\int_{\tau}^{t}\left(\eta_{1}(s)+\frac{\eta_{2}(s)}{1-\widehat{\vartheta }}\right) ds\right) u(\tau).
\label{DESIGUALDAD_GRONWALL_2}
\end{align*}
\end{lemma}
\subsection{Existence and uniqueness for \eqref{SISTEMA_GENERICO_IDEPCAG}}
\begin{theorem}{(\bf{Uniqueness}) (\cite{CTP2019}, Theorem 4.5)}
Consider Lemma \ref{LEMA_GRONWALL_PINTO_2} and the I.V.P for \eqref{sistema_idepcag_general_abstract} with $y(t,\tau,y(\tau))$. Let (H1)-(H2) hold. Then, there exists a unique solution $y$ for \eqref{sistema_idepcag_general_abstract} on $[\tau,\infty)$. Moreover, every solution is stable.
\end{theorem}

\begin{lemma}{(\bf{Existence of solutions in $[\tau,t_{k})$}) (\cite{CTP2019}, Lemma 4.6)}
Consider Lemma \ref{LEMA_GRONWALL_PINTO_2} and the I.V.P for \eqref{sistema_idepcag_general_abstract} with $y(t,\tau,y(\tau))$. Let (H1)-(H2). Then, for each $y_0\in\mathbb{R}^n$ and $\zeta_k\in [t_{k-1},t_{k})$ there exists a solution $y(t)=y(t,\tau,y(\tau))$ of \eqref{sistema_idepcag_general_abstract} on $[\tau,t_r)$ such that $y(\tau)=y_0$.
\end{lemma}

\begin{theorem}{(\bf{Existence of solutions in $[\tau,\infty$}) (\cite{CTP2019}, Theorem 4.7)}
Consider \ref{LEMA_GRONWALL_PINTO_2} and let (H1)-(H2) hold. Then, for each $(\tau,y_0)\in\mathbb{R}_0^+\times\mathbb{R}$, there exists $y(t)=y(t,\tau,y_0)$ for $t\geq \tau$, a unique solution for \eqref{sistema_idepcag_general_abstract} such that $y(\tau)=y_0$.    
\end{theorem}

\section{The nonautonomous homogeneous linear IDEPCAG}
In this section, we will present the nonautonomous homogeneous linear IDEPCAG 
\begin{equation}
\begin{tabular}{ll}
$z^{\prime }(t)=a(t)z(t)+b(t)z(\gamma (t)),$ & $t\neq t_{k}$ \\ 
$\Delta z|_{t=t_{k}}=c_k z(t_{k}^{-}),$ & $t=t_{k}$
\end{tabular}
\label{SISTEMA_IDEPCAG_ORIGINAL_}
\end{equation}
where $z\in\mathbb{R},t\in\mathbb{R},$ $a(t),b(t)$ are real-valued continuous locally integrable functions, $(c_k)_{k\in\mathbb{N}}$ is a real sequence such that $1+c_k\neq 0$ $\forall k\in\mathbb{N}$ and $\gamma (t)$ is a generalized piecewise constant argument.\\

As the oscillatory definition has an asymptotic component, we will only be interested in the forward continuation of solutions. So, during the rest of the work, we will assume $\gamma(\tau):=\tau$ if $t_{k(\tau)}\leq \gamma(\tau)< \tau <t_{k(\tau)+1},$ where $k(\tau)$ is the only $k\in\mathbb{Z}$ such that $t_{k(\tau)}\leq \tau \leq t_{k(\tau)+1},$  and we will adopt the following notation for the matricial products and sums:
\begin{equation*}
\prod_{j=1}^{n}A_j=
    \begin{cases}
        A_n\cdot A_{n-1} \cdots A_1, & \text{if } n\geq 1,\\
        \qquad \quad I & \text{if } n<1.
    \end{cases}
    \quad \text{ and } \quad
\sum_{j=1}^{n}A_j=
    \begin{cases}
        A_1+\ldots +A_{n}, & \text{if } n\geq 1,\\
        \qquad \quad 0 & \text{if } n<1.
    \end{cases}
\end{equation*}

Let $x(t)=\phi(t,\tau)x(\tau)=\displaystyle{\exp\left(\int_{\tau}^{t}a(u)du\right)}x(\tau)$ the  solution of the ordinary differential equation
\begin{eqnarray*}
&&x^{\prime}(t)=a(t)x(t),\\
&&x_0=x(\tau), \qquad t,\tau\in [\tau,\infty).
\end{eqnarray*}
We will assume the following hypothesis:
\begin{enumerate}
    \item[\textbf{(H3)}] Let
$\displaystyle{\rho_{k}^{+}(a) =exp\left(\int_{t_{k}}^{\zeta_{k}}\left| a(u)\right| du\right),}$  $\displaystyle{\rho_{k}^{-}(a)
=exp\left(\int_{\zeta_{k}}^{t_{k+1}}\left| a(u)\right| du\right),}$\\ 
$\rho_{k}(a)=\rho_{k}^{+}(a)\rho_{k}^{-}(a)$ and
$\nu_{k}^{\pm }(b)=\rho_{k}^{\pm }(a)\ln \rho_{k}^{\pm }(b),$

and assume that $$\rho(a)=\displaystyle{\sup_{k\in \mathbb{Z}}\rho_{k}(a)<\infty},\qquad \nu^{\pm}(b)=\sup_{k\in \mathbb{Z}}\nu_{k}^{\pm}(b)<\infty,$$
where
\begin{equation}
\nu_k^{+}(b)<\nu^{+}(b)<1,\quad \nu_k^{-}(b)<\nu^{-}(b)<1. \label{Cond_invertibilidad}
\end{equation}
\end{enumerate}
\noindent Consider the following definitions
\begin{eqnarray}
&j(t,\zeta_{k})&=1+\int_{\zeta_{k}}^{t}\exp{\left(\int_{s}^{\zeta_k}a(u)du\right)b(s)ds},\notag\\   
&e(t,\zeta_k)&=\phi(t,\zeta_{k(t)}) j(t,\zeta_{k(t)})\notag\\
&&=\left(\exp{\int_{\zeta_{k(t)}}^{t}a(u)du}\right)\notag\\
&&\cdot\left(1+\int_{\zeta_{k(t)}}^{t}\exp{\left(\int_{s}^{\zeta_{k(s)}}a(u)du\right)b(s)ds}\right). \label{MATRIZ J}
\end{eqnarray}
\newpage
\begin{remark}
As a consequence of $(H3)$, it is important to notice the following facts: 
\begin{enumerate}
\item[(i)] Due to condition \eqref{Cond_invertibilidad},
$j^{-1}(t_k,\zeta_{k})$ and $j^{-1}(t_{k+1},\zeta_{k})$ are well defined $\forall k\in \mathbb{Z},$ and 
\begin{eqnarray*}
\left|j^{-1}(t_k,\zeta_k)\right| \leq \sum\limits_{k=0}^{\infty } \left[ \nu^{+} (b)\right] ^{k}  =\frac{1}{1-\nu^{+} (b)},\quad 
\left| j(t_{k+1},\zeta_{k})\right| \leq 1+\nu^{-} (b), \label{COTA_J}\\
\left|j^{-1}(t_{k+1},\zeta_k)\right| \leq \sum\limits_{k=0}^{\infty } \left[ \nu^{-} (b)\right] ^{k}  =\frac{1}{1-\nu^{-} (b)},\quad 
\left| j(t_{k},\zeta_{k})\right| \leq 1+\nu^{+} (b).\label{COTA_J2}
\end{eqnarray*}
Also, if we set $t_0=\tau$, we are considering that $j^{-1}(\tau,\gamma(\tau))$ exists.

\item[(ii)] On the other hand, if we want a nonzero solution of a linear IDEPCAG, we need $j(t,\zeta_k)\neq 0, \forall k\in\mathbb{Z}$ and $\forall t\in [\tau,\infty)$ (See Remark \ref{remark_ceros} and \cite{Kuo_pinto_2013}).
\end{enumerate}
\end{remark}

\subsection{The fundamental solution of the homogeneous linear IDEPCAG}
The following results can be found in  \cite{Torres1} and \cite{TORRES_var} in a more general matricial context. They are the IDEPCAG extension of \cite{P2011}. We will use the next scalar version of them:
\begin{theorem}\label{TEO_FORMULA_var_PAram}
Let the following linear IDEPCAG system
\begin{equation}
\begin{tabular}{ll}
$z^{\prime }(t)=a(t)z(t)+b(t)z(\gamma (t)),$ & $t\neq t_{k}$ \\ 
$z(t_{k})=\left( 1+c_k\right) z(t_{k}^{-}),$ & $t=t_{k}$\\
$z_0=z(\tau).$
\end{tabular}
\label{sistema_w}
\end{equation}
If (H3) holds, then the unique solution of \eqref{sistema_w} is
 \begin{equation*}
z(t)=w(t,\tau)z(\tau),  \quad t\in [\tau,\infty),\label{DEFINICION_GENERAL_MATRIZ_W}
\end{equation*}
where $w(t,\tau)$ is given by 
\begin{equation}
w(t,\tau)=w(t,t_{k(t)})\left(\prod_{r=k(\tau)+2}^{k(t)}\left(1+c_{r}\right) w(t_{r},t_{r-1})\right) 
\left(1+c_{k\left(\tau \right) +1}\right) w(t_{k(\tau)+1},\tau )
\label{MATRIZ_FUNDAMENTAL_IDEPCAG}
\end{equation}
for $t\in I_{k(t)}$ and $\tau \in I_{k(\tau) },$
where $w(t,s)$ is defined as 
\begin{equation*}
w(t,s)=e(t,\gamma(s))e^{-1}(s ,\gamma(s)),\qquad \text{if }t,s
\in I_k=[t_{k},t_{k+1}]. \label{MATRIZ_w}
\end{equation*}
In a more explicit form, the unique solution of \eqref{sistema_w} can be written as 
\begin{eqnarray}
z(t) &=&w(t,t_{k(t)})\left(
\prod_{r=k(\tau)+2}^{k(t)}\left(1+c_{r}\right) w(t_{r},t_{r-1})\right)\notag\\
&&\cdot
\left(1+c_{k(\tau)+1}\right) w(t_{k(\tau)+1},\tau)z(\tau)
,\quad \text{for }t\in I_{k(t)},\,\tau \in I_{k(\tau) }, \label{SOLUCION_FINAL_SISTEMA_IDEPCAG_LINEAL}
\end{eqnarray}
and the discrete solution of \eqref{sistema_w} can be also written by 
\begin{equation}
z(t_{k(t)})=\left(\prod_{r=k(\tau)+2}^{k(t)}\left( 1+c_{r}\right) w(t_{r},t_{r-1})\right) \left(1+c_{k(\tau)+1}\right) w(t_{k(\tau)+1},\tau)z(\tau).
\label{SOLUCION_DISCRETA_SISTEMA_Z}
\end{equation}
The expression given by \eqref{MATRIZ_FUNDAMENTAL_IDEPCAG}
is called the fundamental solution for \eqref{sistema_w}.
\end{theorem}
\begin{proof}
    Let $t,\tau \in I_k=[t_{k},t_{k+1})$ for some $k\in \mathbb{Z}.$
In this interval, we are in the presence of the ordinary system
$$z^{\prime }(t)=a(t)z(t)+b(t)z(\zeta_k).$$
So, the unique solution can be written as
\begin{equation}
z(t)=\phi(t,\tau)z(\tau )+\int_{\tau}^{t}\phi(t,s)b(s)z(\zeta_{k})ds.
\label{variacion_parametros_general}
\end{equation}
Keeping in mind $I_k^{+}$, evaluating the last expression at  $t=\zeta_{k}$ we have
\begin{equation*}
z(\zeta_{k})=\phi(\zeta_{k},\tau )z(\tau)+\int_{\tau}^{\zeta_{k}}\phi(\zeta_{k},s)b(s)z(\zeta_{k})ds,
\end{equation*}
Hence, we get
\begin{eqnarray*}
\left(1+\int_{\zeta_{k}}^{\tau}\phi (\zeta_{k},s)b(s)ds\right)z(\zeta_{k}) &=&\phi (\zeta_{k},\tau )z(\tau).
\end{eqnarray*}
I.e
\begin{equation*}
z(\zeta_{k})=j^{-1}(\tau ,\zeta_{k})\phi(\zeta_{k},\tau )z(\tau).
\label{inicial_gamma_general}
\end{equation*}
Then, by the definition of $e(t,\tau)=\phi(t,\tau)j(t,\tau)$, we have
\begin{equation}
z(\zeta_{k})=e^{-1}(\tau ,\zeta_{k})z(\tau).
\label{condicion_inicial_con_gama_y_E}
\end{equation}
Now, from \eqref{variacion_parametros_general} working on $I_k^{-}$, considering $\tau =\zeta_{k}$, we have
\begin{eqnarray*}
z(t)&=&\phi (t,\zeta_{k})z(\zeta_{k})+\int_{\zeta_{k}}^{t}\phi(t,s)b(s)z(\zeta_{k})ds \\
&=&\phi(t,\zeta_{k})\left( 1+\int_{\zeta_{k}}^{t}\phi (\zeta_{k},s)b(s)ds\right) z(\zeta_{k}).
\end{eqnarray*}
I.e., $z(t)=e(t,\zeta_{k})z(\zeta_{k}).$
So, by \eqref{condicion_inicial_con_gama_y_E}, we can rewrite the last equation   
as
\begin{equation}
z(t)=e(t,\zeta_{k})e^{-1}(\tau ,\zeta_{k})z(\tau).
\label{condicion_inicial_sin_w}
\end{equation}

\noindent Then, setting  
\begin{equation*}
w(t,s)=e(t,\gamma(s))e^{-1}(s ,\gamma(s)),\qquad \text{if }t,s
\in I_k=[t_{k},t_{k+1}],  \label{MATRIZ_w_}
\end{equation*}
we have the solution for \eqref{sistema_w} for $t\in I_{k}=[t_{k},t_{k+1}),$
\begin{equation}
z(t)=w(t,\tau)z(\tau). \label{DEFINICION_GENERAL_MATRIZ_W_}
\end{equation}

\noindent Next, if we consider $\tau=t_k, $ and, assuming left side continuity of \eqref{DEFINICION_GENERAL_MATRIZ_W} at $t=t_{k+1}$, we have
\begin{equation*}
z(t_{k+1}^{-})=w(t_{k+1},t_k )z(t_k)
\end{equation*}
Then, applying the impulsive condition to the last equation, we get
\begin{eqnarray}
z(t_{k+1}) &=&\left(1+c_{k+1}\right)w(t_{k+1},t_k)z(t_k).\label{cond_discreta_osc_kuo}
\end{eqnarray}

\noindent The last expression defines a finite-difference equation whose solution is \eqref{SOLUCION_DISCRETA_SISTEMA_Z}. 
Now, by \eqref{DEFINICION_GENERAL_MATRIZ_W_} and the impulsive condition, we have
$$z(t_{k(\tau)+1})=(1+c_{k(\tau)+1})w(t_{k(\tau)+1},\tau)z(\tau).$$
Hence, considering $\tau=t_k$ in \eqref{DEFINICION_GENERAL_MATRIZ_W} and applying \eqref{SOLUCION_DISCRETA_SISTEMA_Z}, 
we get \eqref{SOLUCION_FINAL_SISTEMA_IDEPCAG_LINEAL}.
In this way, we have solved \eqref{sistema_w} on $[\tau,t).$ 
\end{proof}
\newpage
\begin{remark}\label{remark_ceros}
\begin{itemize}
\item[]
\item[(i)] We used the decomposition of $I_k=I_k^+\cup I_k^-$ to define $W$.\\
In fact, we can rewrite \eqref{MATRIZ_FUNDAMENTAL_IDEPCAG} in terms of the advanced and delayed parts using \eqref{MATRIZ_w}:
\begin{eqnarray*}
w(t,\tau)&=&e(t,\zeta_{k(t)})e^{-1}(t_{k(t)},\zeta_{k(t)})\left(\prod_{r=k(\tau)+2}^{k(t)}\left(1+c_{r}\right) e(t_{r},\zeta_{r-1})e^{-1}(t_{r-1},\zeta_{r-1})\right)\\ 
&&\cdot \left(1+c_{k\left(\tau \right) +1}\right) e(t_{k(\tau)+1},\gamma(\tau))e^{-1}(\tau,\gamma(\tau)),\qquad \zeta_r=\gamma(t_r),
\label{MATRIZ_FUNDAMENTAL_IDEPCAG_avance_retardo}
\end{eqnarray*}
for $t\in I_{k(t)}$ and $\tau \in I_{k(\tau) }.$\\

\item[(ii)] From \eqref{SOLUCION_FINAL_SISTEMA_IDEPCAG_LINEAL} and \eqref{SOLUCION_DISCRETA_SISTEMA_Z}, we can write
\begin{eqnarray}
    &z(t)&=w(t,t_{k(t)})z(t_{k(t)}),   \label{z_de_t}
\end{eqnarray}
for $t\in I_k=[t_k,t_{k+1}),$ where $$w(t,t_{k(t)})=\phi(t,t_{k(t)})j(t,\zeta_{k(t)})j^{-1}(t_{k(t)},\zeta_{k(t)}),$$ and  $z(t_{k(t)})$ is the unique solution of the difference equation 
\begin{eqnarray}
    &z(t_{k(t)+1})=&(1+c_{k(t)+1})w(t_{k(t)+1},t_{k(t)})z(t_{k(t)}),\label{z_k_de_t}
\end{eqnarray}
with $z_0=z_{k(\tau)}.$
\item[(iii)]
It is important to notice that a \eqref{sistema_w} has a zero on the interval $I_j=[t_j,t_{j+1})$ at $t=\alpha_j$ if and only if 
\begin{equation*}
1+\displaystyle{\int_{\zeta_{j}}^{\alpha_j}\exp\left(\int_{s}^{\zeta_{j}}a(u)du\right)b(s)ds=0}
\end{equation*}
holds, for some $\alpha_j\in I_j.$
\end{itemize}
\end{remark}


\begin{remark}
\begin{itemize}
\item[]
\item[(i)]   Considering $b(t)=0$, we recover the classical fundamental solution of the impulsive linear differential equation (see \cite{Samoilenko}).
\item[(ii)] If $c_k=0, \forall k\in\mathbb{Z}$, we recover the DEPCAG case studied by M. Pinto in \cite{P2011}. 
\item[(iii)] If we consider $\gamma(t)=p\left[\dfrac{t+l}{p}\right]$ with $p<l$, we recover the IDEPCA case studied by K-S. Chiu in \cite{Kuo2023}.
\end{itemize}
\end{remark}

\section{Oscillations in Linear IDEPCAG}
In this section, we will present sufficient conditions
to ensure the oscillatory character of solutions for the scalar version of IDEPCAG \eqref{sistema_w}.

\subsection{\textbf{Oscillations in scalar Linear IDEPCAG}}

Following the ideas given in \cite{kuo2011}, we will provide useful results concerning oscillatory and nonoscillatory conditions for linear IDEPCAG systems, knowing the oscillatory properties of either the discrete or the continuous solution of a linear IDEPCAG.  
\begin{theorem}
Let the scalar IDEPCAG:
\begin{equation}
\begin{tabular}{ll}
$z^{\prime }(t)=a(t)z(t)+b(t)z(\gamma(t)),$ & $t\neq t_{n}$ \\ 
$z(t_{n})=\left(1+c_n\right) z(t_{n}^{-}),$ & $t=t_{n}$\\
$z_0=z(\tau),$
\end{tabular}
\label{sistema_w_escalar}
\end{equation}
with $1+c_n\neq 0$ for all $n\in\mathbb{N}$.
Then: 
\begin{itemize}
\item[a)] If the discrete solution of \eqref{sistema_w_escalar} (i.e., the solution of the difference equation \eqref{cond_discreta_osc_kuo}) is oscillatory, then the solution $z(t)$ of \eqref{sistema_w_escalar} is oscillatory.
\item[b)] $z(t)$ will be oscillatory if the sequence
\begin{equation*}
\{w_n\}=\left\{(1+c_n)\left(\displaystyle{\dfrac{1+\int_{\zeta_{n}}^{t_{n+1}}\exp\left(\int_{s}^{\zeta_{n}} a(r)dr\right)b(u)du}{1+\int_{\zeta_{n}}^{t_{n}}\exp\left(\int_{s}^{\zeta_{n}} a(r)dr\right)b(u)du}}\right)\right\}_{n\geq n(\tau)+1}
\end{equation*}
is not eventually positive nor eventually negative.
\item[c)] Suppose that the discrete solution of \eqref{sistema_w_escalar} $z(t_n)_{n\geq n(\tau)}$ is nonoscillatory. Then, $z(t)$ is nonoscillatory if and only if 
\begin{equation}
\begin{cases}
1+\displaystyle{\int_{\zeta_{k(t)}}^{t}\exp\left(\int_{s}^{\zeta_{k(t)}}a(u)du\right)b(s)ds}>0, & \label{cond_osci_t}\\
\qquad \qquad \text{ or }\\
1+\displaystyle{\int_{\zeta_{k(t)}}^{t}\exp\left(\int_{s}^{\zeta_{k(t)}}a(u)du\right)b(s)ds}<0 & 
\end{cases}
\end{equation}
hold for $t\in[t_{k(t)},t_{k(t)+1})$, $\forall t,t_k\geq M$, where $k=k(t)$ and $M\in\mathbb{N}$ sufficiently large. 
\end{itemize}
\end{theorem}
\begin{remark}
    Condition \eqref{cond_osci_t} is necessary, because $w(t,t_{k(t)})$ involves
    \begin{eqnarray*}
    &&j(t,\zeta_{k(t)})j^{-1}(t_{k(t)},\zeta_{k(t)})=\left(
    \dfrac{1+\int_{\zeta_{k(t)}}^{t}\exp\left(\int_{s}^{\zeta_{k(t)}} a(r)dr\right)b(s)ds}{1+\int_{\zeta_{k(t)}}^{t_{k(t)}}\exp\left(\int_{s}^{\zeta_{k(t)}} a(r)dr\right)b(s)ds}\right),\\
    &&t,\zeta_{k(t)}\in I_{k}=[t_k,t_{k+1}), k\geq k(\tau),\,\,k\in\mathbb{N}.
    \end{eqnarray*}
\end{remark}
\begin{proof}
\begin{enumerate}
\item[]
\item[]  We will proceed as in \cite{kuo2011} and \cite{Kuo_pinto_2013} with slight changes due to the impulsive effect:\\

\item[a)] If the discrete solution $z(t_{k(t)})$ of \eqref{sistema_w_escalar} is oscillatory, then \eqref{sistema_w_escalar} it will be because the discrete solution is part of the whole solution.\\
Formally, 
by \eqref{SOLUCION_FINAL_SISTEMA_IDEPCAG_LINEAL} with $t=t_{k(t)}$  for $t\in[t_{k(t)},t_{k(t)+1})$, we have
$z(t)=z(t_{k(t)})$. I.e., \eqref{SOLUCION_DISCRETA_SISTEMA_Z}. Then, for the classical definition of difference equations, we have that $z(t)$ is oscillatory because $z(t_{k(t)})$ is oscillatory.\\

\item[b)]  By \eqref{MATRIZ J} and \eqref{SOLUCION_DISCRETA_SISTEMA_Z}, the solution of \eqref{sistema_w_escalar} can be rewritten, in terms of $w_n$ and $w(t_{r+1},t_{r})$ as
\begin{equation*}
z(t_{k(t)})=\exp\left(\int_{t_{k(\tau)+1}}^{t_{k(t)}}a(u)du\right)\left(\prod_{r=k(\tau)+1}^{k(t)}
w_r\right)z(t_{k(\tau)+1}),
\end{equation*}
where $w_r=(1+c_r)w(t_{r+1},t_r)$ and
\begin{eqnarray*}
&w(t_{r+1},t_{r})&=e(t_{r+1},\zeta_{r})e^{-1}(t_{r},\zeta_{r})\\
&&=\exp\left(\int_{t_{r}}^{t_{r+1}}a(u)du\right)
\left(\dfrac{1+\int_{\zeta_{r}}^{t_{r+1}}\exp\left(\int_{u}^{\zeta_{r}} a(r)dr\right)b(u)du}{
1+\int_{\zeta_{r}}^{t_{r}}\exp\left(\int_{u}^{\zeta_{r}} a(r)dr\right)b(u)du}\right).   
\end{eqnarray*}
It is clear that $z(t_{k(t)})$ is oscillatory if $\{w_n\}$ is not eventually positive or negative. Hence, by part $a)$, $z(t)$ is oscillatory if $\{w_n\}_{n\geq n(\tau)+1}$ is not eventually positive nor eventually negative.\\  
\item[c)] First, we do not have to forget that the impulsive effect is already considered in $z(t_n)_{n\geq n(\tau)}$, because we defined the discrete solution in that way.\\

$\left(\Rightarrow\right)$ If $z(t_n)_{n\geq n(\tau)}$ is nonoscillatory, we can suppose $z(t_n)>0$ and $z(t)>0$ for all $t,n\geq M$ with  $M\in\mathbb{N}$ sufficiently large. Also, by \eqref{DEFINICION_GENERAL_MATRIZ_W_} with $\tau=t_n$  for $t\in[t_{n(t)},t_{n(t)+1})$ and $t\geq M$, we have
$$z(t)=w(t,t_n)z(t_n),$$
where
\begin{eqnarray*}
&w(t,t_{n})&=\exp\left(\int_{t_{n}}^t a(u)du\right)
\left(\dfrac{1+\int_{\zeta_{n}}^{t}\exp\left(\int_{s}^{\zeta_{n}} a(u)du\right)b(s)ds}{
1+\int_{\zeta_{n}}^{t_n}\exp\left(\int_{s}^{\zeta_{n}} a(u)du\right)b(s)ds}\right),   
\end{eqnarray*}
for $t\in[t_{n},t_{n+1})$.\\ 
Then, we have $w(t,t_n)>0,$  i.e. \eqref{cond_osci_t}. The case $z(t_n)<0$ and $z(t)<0$ is analog.\\

$\left(\Leftarrow\right)$ Suppose that $z(t_n)>0$, for all $n\geq M$ sufficiently large and \eqref{cond_osci_t} holds. We have to prove that $z(t)$ is nonoscillatory. We will proceed by contradiction: assume that $z(t)$ is oscillatory. Then, according to the definition of oscillatory solutions, there must exist $(a_n),(b_n)$ sequences such that $(a_n),(b_n)\to\infty$ as $t\to \infty$ and $z(a_n)\leq 0 \leq z(b_n)$. Let $t_n < a_n < t_{n+1}$. From \eqref{DEFINICION_GENERAL_MATRIZ_W_} with $\tau=t_n$ and $t=a_n$, we get $$z(a_n)=w(a_n,t_n)z(t_n).$$
Since $z(a_n)\leq 0$  and $z(t_{n})>0$, we conclude that $w(a_n,t_n)\leq 0,$ and this fact contradicts to \eqref{cond_osci_t}. The case $z(t_n)<0$ for $n\geq M$ is analog, taking $z(b_n)\geq 0$ with $t_n < b_n < t_{n+1}.$ I.e.,
$$z(b_n)=w(b_n,t_n)z(t_n).$$
Again, since $z(b_n)\geq 0$  and $z(t_{n})<0$, we conclude that $w(b_n,t_n)\leq 0,$ and this fact also contradicts to \eqref{cond_osci_t}. 
\end{enumerate}
\end{proof}
Next, inspired by \cite{WiAf88}, we will state our main result, a version of the Wiener-Aftabizadeh criterion for oscillatory solutions of \eqref{sistema_w_escalar}:
\begin{theorem}{\textbf{(Main Result: Aftabizadeh-Wiener IDEPCAG 
 oscillatory criterion)}}\label{criterio_idepcag_wiener}
Let the scalar IDEPCAG \eqref{sistema_w_escalar} with $1+c_k\neq 0$ for all $k\geq M$, $M\in\mathbb{N}$ and $M$ sufficiently large. If either of the conditions:    
$$
\text{ if } 1+c_{k(t)}>0 \text{ and }
\begin{cases}
\displaystyle{\lim_{k\to\infty}\sup_{k\in\mathbb{N}}\int_{t_{k(t)}}^{\zeta_{k(t)}}\exp\left(\int_{s}^{\zeta_{k(s)}}a(u)du\right)b(s)ds>1,} & \\
\qquad \qquad \quad \text{ or} \\
\displaystyle{\lim_{k\to\infty}\inf_{k\in\mathbb{N}}\int_{\zeta_{k(t)}}^{t_{k(t)+1}}\exp\left(\int_{s}^{\zeta_{k(s)}}a(u)du\right)b(s)ds<-1,} & 
\end{cases}
$$
or
$$
\text{ if }1+c_{k(t)}<0 \text{ and }
\begin{cases}
\displaystyle{\lim_{k\to\infty}\inf_{k\in\mathbb{N}}\int_{t_{k(t)}}^{\zeta_{k(t)}}\exp\left(\int_{s}^{\zeta_{k(s)}}a(u)du\right)b(s)ds<1,} & \\
\qquad \qquad \quad\text{ or }\\
\displaystyle{\lim_{k\to\infty}\sup_{k\in\mathbb{N}}\int_{\zeta_{k(t)}}^{t_{k(t)+1}}\exp\left(\int_{s}^{\zeta_{k(s)}}a(u)du\right)b(s)ds>-1.} &
\end{cases}
$$
hold for $t\geq k$, with $k=k(t)$ and $\gamma(t)=\zeta_{k(t)}$ if $t\in I_{k(t)}=[t_{k(t)},t_{k(t)+1})$. Then, every solution $z(t)$ of \eqref{sistema_w_escalar} oscillates.
\end{theorem}
\begin{proof}
    If \eqref{sistema_w_escalar} has nonoscillatory solutions, then $z(t)>0$  or $z(t)<0$ for $k=k(t)$ sufficiently large, for all $k(t)\geq M$ sufficiently large and $\gamma(t)=\zeta_{k(t)}$ if $t\in I_{k(t)}=[t_{k(t)},t_{k(t)+1})$. For the proof, we will consider the advanced and delayed intervals $I_k^+=[t_k,\zeta_k)$ and $I_k^-=[\zeta_k,t_{k+1})$, respectively. Let's suppose that $z(t)>0$ for $t\in I_k^-$. Then, it is easy to see that
    \begin{equation}
\left(z(t)\exp\left(\int_t^{\zeta_k}a(u)du\right)\right)^{\prime}=b(t)z(\gamma(t))\exp\left(\int_t^{\zeta_{k}}a(u)du\right). \label{derivada_oscilatorio}
    \end{equation}
Integrating the last expression from $\zeta_{k}$ to $t_{k+1}$ we have
\begin{eqnarray*}
&&
z(t_{k+1}^-)\exp\left(\int_{t_{k+1}}^{\zeta_{k}}a(u)du\right)=z(\zeta_{k})\left(1+\int_{\zeta_{k}}^{t_{k+1}}\exp\left(\int_{s}^{\zeta_{k(s)}}a(u)du\right)b(s)ds\right).
\end{eqnarray*}
Next, applying the impulsive condition $(1+c_{k(t)+1})z(t_{k+1}^-)=z(t_{k+1})$, we get
\begin{eqnarray*}
z(t_{k+1})\exp\left(\int_{t_{k+1}}^{\zeta_k}a(u)du\right)=(1+c_{k+1})z(\zeta_k)\left(1+\int_{\zeta_{k}}^{t_{k+1}}\exp\left(\int_{s}^{\zeta_{k(s)}}a(u)du\right)b(s)ds\right)
\end{eqnarray*}
By hypothesis, $z(t_{k+1}),z(\zeta_k)>0$. Hence, if $(1+c_{k+1})>0$, then we have
\begin{equation*}
    \int_{\zeta_{k}}^{t_{k+1}}\exp\left(\int_{s}^{\zeta_{k(s)}}a(u)du\right)b(s)ds>-1.
\end{equation*}
I.e., 
\begin{equation*}
 \lim_{k\to\infty}\inf_{k\in\mathbb{N}}\int_{\zeta_{k}}^{t_{k+1}}\exp\left(\int_{s}^{\zeta_{k(s)}}a(u)du\right)b(s)ds\geq -1,
\end{equation*}
Analogously, if $1+c_{k+1}<0$, we can conclude
\begin{equation*}
\lim_{k\to\infty}\sup_{k\in\mathbb{N}}\int_{\zeta_{k}}^{t_{k+1}}\exp\left(\int_{s}^{\zeta_{k(s)}}a(u)du\right)b(s)ds\leq -1,
\end{equation*}
which are a contradiction. Hence, the case $t\in I_k^-$ is proved.\\
Next, we proceed with the case $t\in I_k^+=[t_k,\zeta_k]$. In the same way as before, integrating \eqref{derivada_oscilatorio} from $t_{k}$ to $\zeta_{k}$ we have
\begin{eqnarray*}
&&z(t_{k})\exp\left(\int_{t_{k}}^{\zeta_{k}}a(u)du\right)=z(\zeta_{k})\left(1-\int_{t_{k}}^{\zeta_{k}}\exp\left(\int_{s}^{\zeta_{k(s)}}a(u)du\right)b(s)ds\right).
\end{eqnarray*}
Next, applying the impulsive condition $(1+c_{k})z(t_{k}^-)=z(t_{k})$, we get
\begin{eqnarray*}
(1+c_k)z(t_{k}^-)\exp\left(\int_{t_{k}}^{\zeta_k}a(u)du\right)=z(\zeta_k)\left(1-\int_{\zeta_{k}}^{t_{k}}\exp\left(\int_{s}^{\zeta_{k(s)}}a(u)du\right)b(s)ds\right)
\end{eqnarray*}
By hypothesis, $z(t_{k}^-),z(\zeta_k)>0$. Hence, if $(1+c_{k})>0$, then we have
\begin{equation*}
\int_{t_{k}}^{\zeta_{k}}\exp\left(\int_{s}^{\zeta_{k(s)}}a(u)du\right)b(s)ds<1.
\end{equation*}
I.e., 
\begin{equation*} \lim_{k\to\infty}\sup_{k\in\mathbb{N}}\int_{\zeta_{k}}^{t_{k}}\exp\left(\int_{s}^{\zeta_{k(s)}}a(u)du\right)b(s)ds\leq 1,
\end{equation*}
Analogously, if $1+c_{k}<0$, we can conclude
\begin{equation*}
\lim_{k\to\infty}\inf_{k\in\mathbb{N}}\int_{t_{k}}^{\zeta_k}\exp\left(\int_{s}^{\zeta_{k(s)}}a(u)du\right)b(s)ds\geq 1,
\end{equation*}
which are a contradiction. Hence, the case $t\in I_k^+$ is proved.\\
The proof is analog for $z(t) < 0$ and $t$ sufficiently large. Hence,
\eqref{sistema_w_escalar} has oscillatory solutions only. 
\end{proof}
\begin{remark}
\begin{itemize}
\item[]
\item If $c_k=0,\,\forall k\in\mathbb{Z},$ we recover the results of K-S. Chiu and M. Pinto given in \cite{Kuo_pinto_2013}.
 \end{itemize}
 \end{remark}  
\noindent In a very similar way to Theorem \ref{criterio_idepcag_wiener}, we can deduce the following result:
 \begin{corollary}{\textbf{(Aftabizadeh-Wiener IDEPCAG 
 nonoscillatory criterion)}}\label{criterio_idepcag_wiener_nonosc}
Let the scalar IDEPCAG \eqref{sistema_w_escalar} with $1+c_k\neq 0$ for all $k\geq k(t)$. If the conditions:    
$$
\text{ if } 1+c_{k(t)}>0 \text{ and }
\begin{cases}
\displaystyle{\lim_{k\to\infty}\sup_{k\in\mathbb{N}}\int_{t_{k(t)}}^{\zeta_{k(t)}}\exp\left(\int_{s}^{\zeta_{k(s)}}a(u)du\right)b(s)ds\leq 1,} & \\
\qquad \qquad \quad \text{ and} \\
\displaystyle{\lim_{k\to\infty}\inf_{k\in\mathbb{N}}\int_{\zeta_{k(t)}}^{t_{k(t)+1}}\exp\left(\int_{s}^{\zeta_{k(s)}}a(u)du\right)b(s)ds\geq-1,} & 
\end{cases}
$$
or
$$
\text{ if }1+c_{k(t)}<0 \text{ and }
\begin{cases}
\displaystyle{\lim_{k\to\infty}\inf_{k\in\mathbb{N}}\int_{t_{k(t)}}^{\zeta_{k(t)}}\exp\left(\int_{s}^{\zeta_{k(s)}}a(u)du\right)b(s)ds\geq 1,} & \\
\qquad \qquad \quad\text{ and }\\
\displaystyle{\lim_{k\to\infty}\sup_{k\in\mathbb{N}}\int_{\zeta_{k(t)}}^{t_{k(t)+1}}\exp\left(\int_{s}^{\zeta_{k(s)}}a(u)du\right)b(s)ds\leq -1.} &
\end{cases}
$$
hold for $t\geq k$, with $k=k(t)$ sufficiently large and $\gamma(t)=\zeta_{k(t)}$ if $t\in I_{k(t)}=[t_{k(t)},t_{k(t)+1})$, then every solution $z(t)$ of \eqref{sistema_w_escalar} is nonscillatory.
 \end{corollary}

\section{Some Examples of Oscillatory Linear IDEPCAG systems}
\begin{example}
Let the following IDEPCA: 
\begin{align}
x'(t)&=(\alpha-1)x([t]),\qquad t\neq n \label{ejemplo_intro} \\
x(n)&=\beta x(n^{-}),\qquad t=n, \quad n\in\mathbb{N},
\end{align}
where $t_k=\zeta_k=k,\,\,k\in\mathbb{N}$, and $\alpha,\beta\neq 1$ (the interesting cases).
In virtue of Theorem \ref{criterio_idepcag_wiener}, the solutions of \eqref{ejemplo_intro} are oscillatory if 
$\alpha>0$ and $\beta<0$ or
$\alpha<0$ and $\beta>0.$
I.e. $\alpha \beta<0.$
\end{example}
\begin{example}\label{ejemplo_ladas}
In \cite{Gyori_Ladas}, I. Gyori and G. Ladas studied the following DEPCAs
\begin{eqnarray}
    &x'(t)+px(t)+q_0x([t])=0,\quad p,q_0\in\mathbb{R},  \label{GL1}\\
    &y'(t)+py(t)+q_{-1}y([t-1])=0, \quad p,q\in\mathbb{R}, \label{GL2}
\end{eqnarray}
and they gave the following theorems, using some tools as the characteristic equation for DDE:
\begin{corollary}{(\cite{Gyori_Ladas}(Cor. 8.1.1 and 8.1.2))}
    \begin{enumerate}
        \item Every solution of \eqref{GL1} oscillates if and only if $q_0\geq\dfrac{p}{e^p-1}.$
        \item Every solution of \eqref{GL2} oscillates if and only if $q_{-1}>\dfrac{pe^{-p}}{4(e^p-1)}.$
    \end{enumerate}
\end{corollary}
For \eqref{GL1}, applying Theorem \ref{criterio_idepcag_wiener}  with $\gamma(t)=[t]$, $a(t)=-p,$ $b(t)=-q_0$ and $c_k=0,$ $\forall k\in\mathbb{Z},$  we obtain the following estimation:
\begin{equation*}  \dfrac{q_0}{p}\left(e^{p((k+1)-k)}-1\right)>1\Rightarrow  q_0>\dfrac{p}{e^p-1}.
\end{equation*}
Finally, for \eqref{GL2}, also using Theorem \ref{criterio_idepcag_wiener} and  
\begin{eqnarray*}
&&\displaystyle{\int_{k-1}^{k+1} e^{-p(\gamma(k(s))-s))}ds}= \displaystyle{\int_{k-1}^{k}e^{-p((k-2)-s))}ds+\int_{k}^{k+1}e^{-p((k-1)-s))}ds},
\end{eqnarray*}
we have $
q_{-1}>\dfrac{pe^{-p}}{2(e^p-1)}.$
The bonds found satisfy the ones given in \cite{Gyori_Ladas}.
\end{example}

\begin{remark}\label{remark_ladas}
\begin{itemize}
\item[] 
\item The estimations obtained by I.Gyori and G. Ladas for \eqref{GL1} and \eqref{GL2} are sharper than ours. Still, we recover the same results without analyzing the roots of a specific characteristic equation.
\item In the light of the previous results, it seems that any $\gamma(t)$ such that $\gamma(t)=\zeta_k\notin I_k=[t_k,t_{k+1})$ for $t\in I_k$ could satisfy Theorem \ref{criterio_idepcag_wiener}. ($\gamma(t)=[t-1]$ implies $\zeta_k=k-1$ if $t\in [k,k+1), \, k\in\mathbb{Z}$). This fact deserves future research.
\end{itemize}
\end{remark}
\begin{example}
Let the following scalar linear IDEPCAG
\begin{equation}
    \begin{tabular}{ll}
$z'(t)=a(t)\left(z(t)-z(\gamma(t))\right),$ & $t\neq k$ \\ 
$z(k)=c_k z(k^{-}),$ & $t=k, \quad k\in\mathbb{Z},$ 
\end{tabular}
\label{berek_sistema_generalizado}
\end{equation}
where $a(t)$ is a continuous locally integrable function, $c_k\neq 0,\,\forall k\in\mathbb{Z}$,  $\gamma(t)$ is any generalized piecewise constant argument, i.e., $\gamma(t)=\zeta_k,$ if $t\in I_k=[t_k,t_{k+1}),\,k\in\mathbb{Z}$ and $t_k\leq \zeta_k\leq t_{k+1}$.

As $$\displaystyle{\left(\exp\left(-\int_{\zeta_k}^{s}a(u)du\right)\right)^{\prime}=\exp\left(-\int_{\zeta_k}^{s}a(u)du\right)(-a(s)),}$$
by Theorem \eqref{criterio_idepcag_wiener},
we have the only possible cases occur when $c_{k(t)}<0$ and they are when $\displaystyle{\lim_{k\to\infty}\inf_{k\in\mathbb{N}}\left(1-\exp\left(-\int_{\zeta_k}^{t_{k}}a(u)du\right)\right)<1, }$ or $\displaystyle{\lim_{k\to\infty}\sup_{k\in\mathbb{N}}\left(1-\exp\left(-\int_{t_{k+1}}^{\zeta_k}a(u)du\right)\right)<1.}$  
Hence, if one of the last conditions holds, then the solutions of \eqref{berek_sistema_avanzado} oscillate.\\

\noindent  As we will see, \eqref{berek_sistema_generalizado} is easily solvable.
Let $\phi(t) $ be the fundamental solution of 
\begin{equation}
x^{\prime }(t)=a(t)x(t).  \label{sistema_ordinario_berek}
\end{equation}
It is well known that $\phi ^{-1}(t) $ is the fundamental solution of the adjoint system associated with \eqref{sistema_ordinario_berek}. So, it satisfies 
\begin{equation*}
\left(\phi ^{-1}\right) ^{^{\prime }}(t)=-\phi ^{-1}(t)(a(t)).
\end{equation*}
Therefore, we have
\begin{align*}
j(t,t_k)& =1-\int_{t_k}^{t}\phi(t_k,s)a(s)ds \\
& =1+\phi(t_k)\left(\int_{t_k}^{t}-\phi^{-1}(s)a(s)ds\right) \\
& =1+\phi(t_k)\left(\phi^{-1}(t)-\phi^{-1}(t_k)\right) \\
& =\phi^{-1}(t,t_k).
\end{align*}
As $e(t,t_k)=\phi(t,t_k)j(t,t_k)$, we have $e(t,t_k)=1$ (see Theorem \ref{TEO_FORMULA_var_PAram}).\\
Hence, the solution of \eqref{berek_sistema_generalizado} is 
$$z(t)=\left(\prod_{j=k(\tau)+1}^{k(t)}c_j\right)z(\tau),\quad t\geq \tau.$$
(See example \ref{ejemplo_berek_avanzado}).
The last conclusion applies no matter what locally integrable function $a(t)$ is used.
\begin{figure}[h!]
\centering
\includegraphics[scale=0.25]{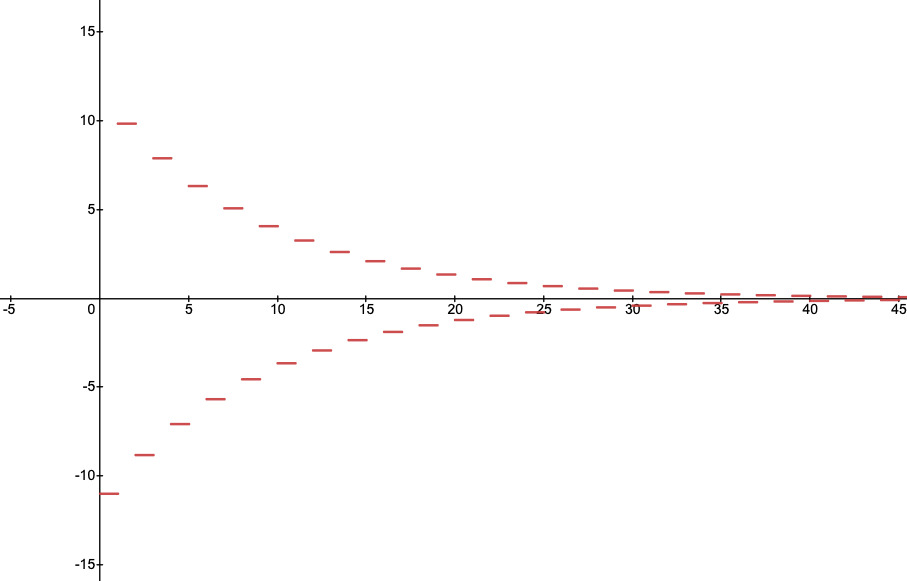}
\caption{Solution of \eqref{berek_sistema_generalizado} with $c_k=-\frac{60}{67}$ and $z(0)=-11.$}
\end{figure}
\begin{remark}
The  IDEPCA presented in \eqref{berek_sistema_avanzado} is a particular case of the last IDEPCAG, with $\gamma(t)=[t+1]$.
\end{remark}
\end{example}

\begin{example}
Let the following scalar linear IDEPCAG
\begin{equation}
    \begin{tabular}{ll}
$x'(t)=-ax(t)+sen(2\pi t)x([t]),$ & $t\neq k,$ \\ 
$x(k)=Cx(k^{-}),$ & $t=k,\quad k\in\mathbb{N},$\\
$x(0)=1.$
\end{tabular}
\label{ejemplo_cos}
\end{equation}
We see that $\displaystyle{\int_{k}^{k+1}\exp\left(-a(k-s)\right)sin(2\pi s)\,\,ds=\dfrac{2\pi}{a^2+4\pi^2}(1-e^a).}$\\
By Theorem \ref{criterio_idepcag_wiener} and Theorem \ref{criterio_idepcag_wiener_nonosc} we have
$$
\begin{cases}
\text{ if }C>0 \text{ and }  
\begin{cases}
a>2.07553 \Rightarrow 
\text{ All the solutions are oscillatory.}\\
a\leq 2.07553 \Rightarrow 
\text{ All the solutions are nonoscillatory.}\\
\end{cases}\\
\text{ if }C<0, 
\text{ all the solutions are oscillatory.}\\
\end{cases}
$$ 
\begin{figure}[h!]
\centering
\includegraphics[scale=0.3]{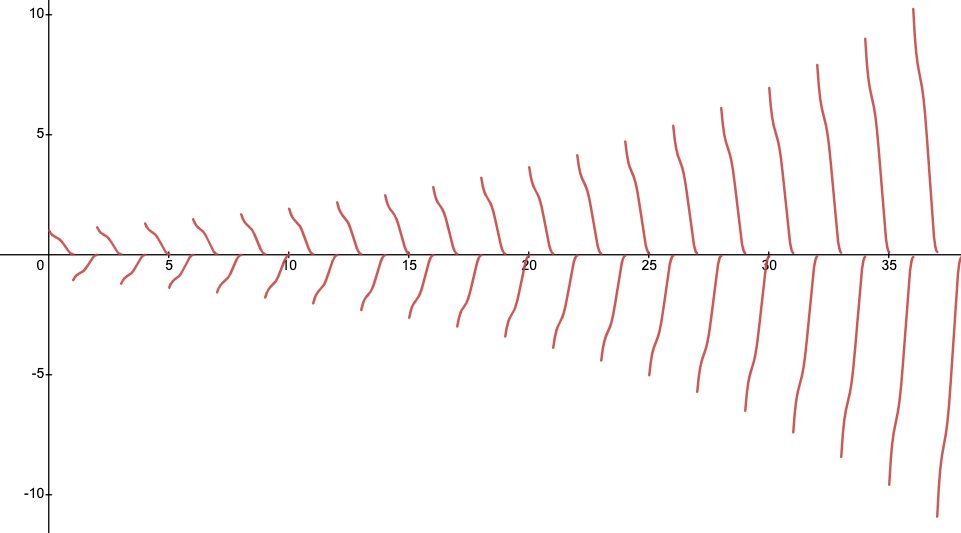}
\caption{Solution of \eqref{ejemplo_cos} with $x(0)=1, C=-100$ and $a=1.998$ .}
\end{figure}
\end{example}
\section*{Acknowledgements}
The author sincerely thanks Prof. Manuel Pinto for encouraging me to work in this subject and for all his support during my career.
The author also thanks DESMOS PBC for granting permission to use the images employed in this work. They were created with the DESMOS graphic calculator \\ \url{https://www.desmos.com/calculator}.

\end{document}